\documentclass[twoside,10pt]{article}

\usepackage{psfrag}
\usepackage{graphicx}
\usepackage{amssymb}
\usepackage{latexsym}
\usepackage{amsmath}

\newtheorem{theorem}{Theorem}
\newtheorem{proposition}{Proposition}
\newtheorem{definition}{Definition}
\newtheorem{lemma}{Lemma}

\newtheorem{example}{Example}

\date{}

\begin{document}

\title{Isometries of spaces of convex compact subsets of $CAT(0)$-spaces}
\maketitle

\begin{center}
{\large Thomas Foertsch *\footnote{* Supported by the Deutsche Forschungsgemeinschaft (FO 353/1-1)}} 
 \\
\end{center}

\begin{abstract}
In the present paper we characterize the surjective isometries of the space of compact, convex subsets of
proper, geodesically complete $CAT(0)$ spaces in which geodesics do not split, endowed with the Hausdorff metric.
Moreover, an analogue characterization of the surjective isometries of the space of compact subsets of
a proper, uniquely geodesic, geodesically complete metric space in which geodesics do not split, when endowed with the Hausdorff metric, 
is given.
\end{abstract}


\section{Introduction}

Let $(X,d)$ be a metric space. For $A\subset X$, $r>0$ we define the closed tubular neighborhood
$N_r(A)$ of $A$ of radius $r$ as
\begin{displaymath}
N_r(A) \; := \; \Big\{ x\in X \; \Big| \; \exists a\in A \; \mbox{with} \; d(a,x)\le r \Big\} .
\end{displaymath}
For $p\in X$ we also write $B_r(p):=N_r(p)$. The sphere $S_r(p)$ of radius $r$ around $p$ is defined via 
$S_r(p):=\{ x\in X | d(x,p)=r\}$. \\
On the set $\mathcal{B}=\mathcal{B}(X,d)$ of closed, bounded subsets of $X$ the map $d_{H}:X\times X\longrightarrow \mathbb{R}_0^+$ given via
\begin{equation} \label{eqn-def-dH}
d_H (A,B) \; := \; \inf \Big\{ r \; \Big| \; A\subset N_r(B) \;\; \wedge \;\; B\subset N_r(A)\Big\}
\end{equation}
defines the so called Hausdorff metric on $\mathcal{B}(X,d)$. \\

In the late 70's and early 80's several authors started to investigate the relations of isometries of the Euclidean 
space $\mathbb{E}^n$ and those of the space $\mathcal{C}(\mathbb{E}^n)$ of its compact, convex subsets when endowed with the 
Hausdorff metric. Of course, given an isometry $i$ of the Euclidean space, one derives an isometry $I$ of the space
$(\mathcal{C}(\mathbb{E}^n),d_H)$ by setting
\begin{displaymath}
I(C) \; := \; i(C) \hspace{1cm} \forall C\in \mathcal{C}(\mathbb{E}^n).
\end{displaymath} 

In \cite{sch} the author showed that these are the only surjective isometries of $(\mathcal{C}(\mathbb{E}^n),d_H)$. In \cite{g}
it was shown that the same holds for the surjective isometries of $(\mathfrak{C}(\mathbb{E}^n),d_H)$,
where $\mathfrak{C}(\mathbb{E}^n)$ denotes the set of compact subsets of $\mathbb{E}^n$, and \cite{gt} generalizes
these observations to certain non-Euclidean cases. Here the authors raise the question whether a similar statement also holds for 
real hyperbolic spaces and we are not aware of the fact that this has been considered anywhere so far. \\

The main purpose of this paper are the following broad generalizations of Schmidt's Theorem for
the space of compact convex subsets of the Euclidean space and Gruber's Theorem for the space of
compact subsets of the Euclidean space:

\begin{theorem} \label{theo-main}
Let $(X,d)$  be a proper, uniquely geodesic, geodesically complete metric space such that geodesics do not split and assume
that the unique midpoint map $m$ of $(X,d)$ is convex. Let further $I$ be a surjective isometry of 
$(\mathcal{C}(X,d),d_H)$. Then there exists an isometry $i\in Isom(X,d)$ such that
\begin{displaymath}
I(C) \; = \; i(C) \hspace{1cm} \forall C\in \mathcal{C}(X,d) .
\end{displaymath}
\end{theorem}


\begin{theorem} \label{theo-main-2}
Let $(X,d)$ be a proper, uniquely geodesic, geodesically complete metric space in which geodesics do not split and $I$
be an isometry of $(\mathfrak{C}(X,d),d_H)$ onto itelf. Then there exixts an isometry $i\in Isom(X,d)$ such that 
\begin{displaymath}
I(C) \; = \; i(C) \hspace{1cm} \forall C\in \mathfrak{C}(X,d) .
\end{displaymath}
\end{theorem}

For the precise definitions of those properties characterizing the metric spaces as considered in Theorem \ref{theo-main}
and Theorem \ref{theo-main-2} we refer the reader to the Sections \ref{subsec-basics} and \ref{subsec-convex-midpoint}. 
Note, however, that our Theorems in particular apply to all 
proper, geodesically complete $CAT(0)$-spaces in which geodesics do not spit, therefore for instance to  all complete,
connected, simply connected Riemannian manifolds of non-positive curvature and, moreover, to all finite dimensional Banach spaces with
strictly convex norm balls.


{\bf Outline of the paper:} In Section \ref{subsec-basics} we recall some definitions and set up the notation we are frequently going to use
in this paper. In Section \ref{subsec-convex-midpoint} convex midpoint maps are introduced examples of which will be given
in Section \ref{subsec-closed-convex-sets}, where we also observe that one consequence of our Theorem \ref{theo-main} 
is the existence of a certain class of geodesics, which is invariant under isometries of the spaces considered. 
Here we also point out that this can be interpreted as a Mazur-Ulam type theorem for metric spaces. \\
Then, in Section \ref{sec-proof-1}, we prove Theorem \ref{theo-main}, while the proof of Theorem \ref{theo-main-2} is subject to
Section \ref{sec-proof-2}. \\

{\bf Acknoledgements:} It is a pleasure to thank Anders Karlsson, Viktor Schroeder and, in particular, Mario Bonk for usefull discussions.


\section{Preleminaries}

\subsection{Basic definitions and notation}
\label{subsec-basics}

Recall that a metric space is called {\it proper} if and only if all its closed metric balls are compact. Further recall that for
$a,b\in \mathbb{R}$, $I\in \{ (-\infty ,b],[a,b],[a,\infty )\}$ an isometric embedding $\gamma :(I,|\cdot |)\longrightarrow (X,d)$
is called a geodesic of $(X,d)$. In case $I=[a,b]$ we say that $\gamma$ connects $\gamma (a)$ to $\gamma (b)$, while
for $I=[a,\infty )$ we also refer to $\gamma$ as a geodesic ray initiating in ${\gamma}(a)$. \\
A metric space $(X,d)$ is said to be {\it geodesic}, if and only if for each $x,y\in X$ there exists a geodesic of $(X,d)$ connecting
$x$ to $y$. Such a geodesic connecting $x$ to $y$ will be denoted by ${\gamma}_{xy}$. Note, however, that in a general geodesic
metric space such a geodesic might not be unique. \\
We call a geodesic metric space $(X,d)$ {\it geodesically complete}, if and only if for each geodesic ${\gamma}_{xy}$ connecting $x\in X$
to $y\in X$ there exists a biinfinite extension, i.e. a geodesic ${\gamma} :(-\infty ,\infty )\longrightarrow X$ such that 
$im \{ {\gamma}_{xy}\} \subset im \{ \gamma \}$. If for each geodesic in $(X,d)$ the image of this biinfinite extension is unique, 
we say that {\it geodesics do not split}. \\
A subset $C\subset X$ of a metric space $(X,d)$ is called {\it convex}, if and only if with two points $a,b\in C$ it also contains the
images of all geodesics connecting $a$ to $b$. By $\mathcal{C}(X,d)$ we denote the set of convex, compact subsets of $(X,d)$,
$\mathfrak{C}(X,d)$ denotes the set of compact subsets of $(X,d)$ and
$Isom (X,d)$ is the group of isometries of $(X,d)$ onto itself. \\
Note that for $p\in X$ we have $\{ p\}\in \mathcal{C}(X,d),\mathfrak{C}(X,d)$ and by a slight abuse of notation we will
also write $p=\{ p\}$. \\
Finally recall that a {\it $CAT(0)$-space} is a geodesic metric space $(X,d)$ such that each points 
$a\in im\{ {\gamma}_{xy}\}$ and $b\in im\{ {\gamma}_{xz}\}$ on a geodesic triangle $\Delta ({\gamma}_{xy},{\gamma}_{xz},{\gamma}_{yz})$
with vertices $x,y,z\in X$ lie not further apart than their corresponding comparison points $\bar{a},\bar{b}\in \mathbb{E}^2$ in 
a comparison triangle $\Delta ({\gamma}_{\bar{x}\bar{y}},{\gamma}_{\bar{x}\bar{z}},{\gamma}_{\bar{y}\bar{z}})$ in $\mathbb{E}^2$.
Here a comparison triangle for $\Delta ({\gamma}_{xy},{\gamma}_{xz},{\gamma}_{yz})$ is a geodesic triangle in $\mathbb{E}^2=(\mathbb{R}^2,d_e)$
with vertices $\bar{x},\bar{y},\bar{z}\in \mathbb{E}^2$ such that 
$d(x,y)=d_e(\bar{x},\bar{y})$, $d(x,z)=d_e(\bar{x},\bar{z})$ and $d(y,z)=d_e(\bar{y},\bar{z})$, and the 
comparison points $\bar{a}\in im \{ {\gamma}_{\bar{x}\bar{y}}\}$ and $\bar{b}\in im \{ {\gamma}_{\bar{x}\bar{z}}\}$ are determined
via $d(x,a)=d_e(\bar{x},\bar{a})$ and $d(x,b)=d_e(\bar{x},\bar{b})$. \\

Note that in a $CAT(0)$ space geodesics connecting two points indeed are unique and metric balls are convex.


\subsection{Convex midpoint maps}

\label{subsec-convex-midpoint}

In this section we introduce the notion of (convex) midpoint maps in a metric space $(X,d)$. Assuming that the underlying metric space
$(X,d)$ is complete, such a midpoint map corresponds to a certain class of geodesics in $(X,d)$.

\begin{definition} \label{def-conv-midp}
Let $(X,d)$ be a metric space. A symmetric map $m:X\times X\longrightarrow X$ is called a midpoint map for $(X,d)$ if
\begin{displaymath}
d\Big( m(x,y),x\Big) \; = \; \frac{1}{2} d(x,y) \; = \; d\Big( m(x,y),y\Big) \hspace{0.5cm} \forall x,y\in X.
\end{displaymath}
Furthermore, the midpoint map $m$ is called convex if
\begin{displaymath}
d\Big( m(x_1,y_1),m(x_2,y_2)\Big) \; \le \; \frac{1}{2} \Big[ d(x_1,x_2)+d(y_1,y_2)\Big] \hspace{0.8cm} \forall x_1,x_2,y_1,y_2\in X.
\end{displaymath}
\end{definition}

\begin{definition} \label{def-m-dist-conv}
Let $(X,d)$ be a metric space and $m:X\times X\longrightarrow X$ be a midpoint map for $(X,d)$. Then $(X,d)$ is said to be
\begin{description}
\item[(i)] $m$-distance convex if 
\begin{displaymath}
d(m(x,y),z) \; \le \; \frac{1}{2} \Big[ d(x,z)+d(y,z)\Big] \hspace{1cm} \forall x,y,z\in X.
\end{displaymath}
\item[(ii)] $m$-global non positively Busemann curved ($m$-global NPBC) if
\begin{displaymath}
d\Big( m(z,x),m(z,y)\Big) \; \le \; \frac{1}{2} d(x,y) \hspace{1cm} \forall x,y,z\in X.
\end{displaymath}
\end{description}
\end{definition}
For an investigation of the notion of distance convexity we refer the reader to \cite{f}. \\

The following lemma is a simple consequence of the Definitions \ref{def-conv-midp} and \ref{def-m-dist-conv}:
\begin{lemma} \label{lemma-necessary}
Let $(X,d)$ be a metric space and $m:X\times X\longrightarrow X$ be a convex midpoint map. Then
\begin{description}
\item[(1)] $m$ is continuous,
\item[(2)] $(X,d)$ is $m$-distance convex and
\item[(3)] $(X,d)$ is $m$-global NPBC.
\end{description}
\end{lemma}

In fact $(X,d)$ being $m$-global NPBC is a sufficient condition for the midpoint map $m$ to be convex:
\begin{lemma}
Let $(X,d)$ be a metric space and $m:X\times X\longrightarrow X$ be a midpoint map for $(X,d)$. Then $m$ is
a convex midpoint map if and only if $(X,d)$ is $m$-global NPBC.
\end{lemma} 
{\bf Proof:} Due to Lemma \ref{lemma-necessary} we only have to show that $(X,d)$ being $m$-global NPBC
is a sufficient condition for $m$ being convex. Let therefore $x_1,x_2,y_1,y_2\in X$, then one has
\begin{eqnarray*}
d\Big( m(x_1,y_1),m(x_2,y_2)\Big) & \le & d\Big( m(x_1,y_1), m(x_2,y_1)\Big)  \; + \; d\Big( m(x_2,y_1),m(x_2,y_2) \Big) \\
& \le & \frac{1}{2} \Big[ d(x_1,x_2) + d(y_1,y_2)\Big] .
\end{eqnarray*}
Thus $m$ indeed is convex. \hfill $\Box$ \\

\begin{example} \label{example-midpoint}
Let $(V,||\cdot ||)$ be a normed vector space. Then
\begin{displaymath}
\mathfrak{m}(x,y) \; := \; \frac{x+y}{2} \hspace{1cm} \forall x,y\in V
\end{displaymath}
is a convex midpoint map.
\end{example}
If $V$ is finite dimensional, then it is not hard to see that $\mathfrak{m}$ as defined above is the only convex 
midpoint map in $(V,||\cdot ||)$. Whether or whether not this generalizes to infinite dimensions is not known to 
the author. \\
Given a convex midpoint map $m$ in a metric space $(X,d)$ and an isometry $I\in Isom(X,d)$, obviously $I\circ m$
again is a convex midpoint map. Thus, establishing the uniqueness of a convex midpoint map in a complete metric space
$(X,d)$ gives rise to a class of distinguished geodesics which is invariant under any isometry $I\in Isom(X,d)$. \\
Unfortunately the author is not aware of a metric space admitting two different convex midpoint maps. However, in case
two such midpoint maps exist in a metric space, then there are infinitely many:

\begin{lemma}
Let $(X,d)$ be a metric space and $m_1,m_2:X\times X\longrightarrow X$ be two convex midpoint maps for $(X,d)$. Then the
map $\tilde{m}:X\times X\longrightarrow X$ defined via
\begin{displaymath}
\tilde{m}(x,y) \; := \; m_1\Big( m_1(x,y),m_2(x,y)\Big) \hspace{1cm} \forall x,y\in X
\end{displaymath}
also is a convex midpoint map for $(X,d)$.
\end{lemma}
{\bf Proof:} That $\tilde{m}$ is a midpoint map simply follows from the $m_1$-distance convexity of $(X,d)$ 
(Lemma \ref{lemma-necessary}). The convexity of $\tilde{m}$ follows from
\begin{eqnarray*}
& & d\Big( \tilde{m}(x_1,y_1),\tilde{m}(x_2,y_2)\Big) \\
& \le & \frac{1}{2}\Big[ d\Big( m_1(x_1,y_1),m_1(x_2,y_2)\Big) \; + \; d\Big( m_2(x_1,y_1),m_2(x_2,y_2)\Big) \Big] \\
& \le & \Big[ d(x_1,x_2)+d(y_1,y_2)\Big] \hspace{1cm} \forall x_1,x_2,y_1,y_2\in X.
\end{eqnarray*}
\hfill $\Box$


\subsection{Spaces of closed, bounded, convex sets}

\label{subsec-closed-convex-sets}

Note that for a proper metric space we have $\mathcal{B}=\mathfrak{C}(X,d)$.  

Given a midpoint map $m$ for $(X,d)$, we call a set $A\subset X$ {\it $m$-convex}, if with two points $a,a'\in A$ it also contains
their $m$-midpoint: $m(a,a')\in A$. The {\it $m$-convex hull} $conv_m(A)$ of a set $A\subset X$ is defined via
\begin{displaymath}
conv_m(A) \; := \; \cap \{ C \; | \; C \; \mbox{is closed and} \; m\mbox{-convex} \; \wedge \; A\subset C\} .
\end{displaymath} 
Denoting by ${\mathcal C}_m$ the subset of ${\mathcal B}$, the elements of which are $m$-convex, we write 
$conv_m :{\mathcal B}\longrightarrow {\mathcal C}_m$ for the map which associates to an $A\in {\mathcal B}$ its 
$m$-convex hull $conv_m(A)$. A function $f:X\longrightarrow \mathbb{R}$ is called {\it $m$-convex} if
\begin{displaymath}
f\Big( m(x,y)\Big) \; \le \; \frac{1}{2} \Big[ f(x) \; + \; f(y)\Big] \hspace{1cm} \forall x,y\in X.
\end{displaymath} 

With this terminology it is easy to prove the 
\begin{lemma} \label{lemma-m-conv-dist}
Let $(X,d)$ be a metric space and $m:X\times X\longrightarrow X$ be a convex midpoint map for $(X,d)$.
Let further $C\subset X$ be a closed $m$-convex set in $(X,d)$. Then the map 
\begin{displaymath}
d_C : X\longrightarrow \mathbb{R}^+_0 , \hspace{1cm} x\longmapsto dist(x,C)
\end{displaymath}
is $m$-convex.
\end{lemma}

Furthermore one obtains the
\begin{lemma} \label{lemma-bs} \cite{bs}
The map $conv_m : {\mathcal B} \longrightarrow {\mathcal C}$ is $1$-Lipschitz and does not change the diameter. 
\end{lemma} 
{\bf Proof:} Connecting $b,b'\in B\in {\mathcal B}$ by the $m$-geodesic segment increases neither $diam B$ nor
the Hausdorff distance to any $B'\in {\mathcal B}$ by convexity of the distance function. The claim follows
since $conv_m(B)$ coincides with the closure of $\lim_{n\rightarrow \infty}{\cup}_n B_n $, where $B_0:=B$
and $B_{n+1}$ is obtained from $B_n$ by connecting each pair of points $b,b'\in B_n$ by the $m$-geodesic segment.
\hfill $\Box$

\begin{proposition} \label{prop-midpoint-1}
Let $(X,d)$ be a metric space and $m:X\times X\longrightarrow X$ be a convex midpoint map for $(X,d)$.
Then the map $M:{\mathcal C}_m\times {\mathcal C}_m\longrightarrow {\mathcal C}_m$ defined via
\begin{displaymath}
M(A,A') \; := \; conv_m \Big( \Big\{ x\in X \; \Big| \; \exists a\in A, a'\in A' \; \mbox{such that} \; x=m(a,a')\Big\} \Big)
\;\;\;\;\; \forall A,A' \in {\mathcal C}_m
\end{displaymath} 
is a convex midpoint map for $({\mathcal C}_m,d_H)$.
\end{proposition}
{\bf Proof:} (1) {\it $M$ is a midpoint map}:  Let $\tilde{M}\subset X$ be the set of the midpoints $m(a,a')$
for all $a\in A$, $a'\in A'$. We set $\lambda := \frac{1}{2}d_H(A,A')$ and assume that there exists $b\in \tilde{M}$
with $dist(b,A)>\lambda$. This $b$ is a midpoint $b=m(a,a')$ with $a\in A$ and $a'\in A'$. Since $A$ is $m$-convex, 
the distance function $d_A$ is $m$-convex (Lemma \ref{lemma-m-conv-dist}). Thus $dist(a',A)\ge 2dist(b,A)$, because 
$dist(a,A)=0$. Hence, $dist(a',A)>d_H(A,A')$ contradicting the definition of $d_H(A,A')$. This shows that $\tilde{M}$
lies in the closed $\lambda$-neighborhood of $A$, $N_{\lambda}(A)$. \\
On the other hand, for each $a\in A$ there is $b\in \tilde{M}$ with $d(b,a)\le \lambda$: let $b:=m(a,a')$, where
$a'\in A'$ is the closest point to $a$, thus $d(a,a')\le 2\lambda$. This shows that
$A\subset N_{\lambda}(\tilde{M})$. Thus $d_H(\tilde{M},A)\le \lambda$ and, similarly, $d_H(\tilde{M},A')\le \lambda$. 
By the triangle inequality we have $2\lambda \le d_H(A,\tilde{M})+d_H(\tilde{M},A')\le 2\lambda$ and hence
\begin{displaymath}
d_H(\tilde{M},A) \; = \; \lambda \; = \; d_H(\tilde{M},A').
\end{displaymath}   
For the $m$-convex hull $M=conv_m(\tilde{M})$ we have $d_H(\tilde{M},A),d_H(\tilde{M},A')\le \lambda$ by Lemma \ref{lemma-bs}.
Hence, $d_H(M,A)=\lambda =d_H(M,A')$ and $M$ indeed is a midpoint map for $(C_m,d_H)$. \\
(2) {\it $M$ is convex}: We need to show that
\begin{displaymath}
d_H\Big( M(A,A'),M(B,B')\Big) \; \le \; \frac{1}{2} \Big[ d_H(A,B) + d_H(A',B')\Big] \hspace{0.7cm} 
\forall A,B,A',B'\in \mathcal{C}_m .
\end{displaymath}
Therefore let $A,A',B,B'\in \mathcal{C}_m$ and set $r_0:=d_H(A,B)$ and $r_0':=d_H(A',B')$. All we have to prove now is, 
that given an arbitrary $x\in M(A,A')$, there exists a $y\in M(B,B')$ such that $d(x,y)\le \frac{r_0+r_0'}{2}$. \\
(i) Suppose first that to $x\in M(A,A')$ there exist $a\in A$ and $a'\in A'$ such that $x=m(a,a')$. Due
to the definition of $r_0$ and $r_0'$ there exist $b\in B$ and $b'\in B'$ such that $d(a,b)\le r_0$ and
$d(a',b')\le r_0'$. Now $m$ is a convex midpoint map for $(X,d)$ and for $y:=m(b,b')\in M(B,B')$ we find
\begin{eqnarray*}
d(y,x) & = & d\Big( m(b,b'), m(a,a')\Big) \\
& \le & \frac{1}{2} \Big[ d(b,a)+d(b',a')\Big] \\
& \le & \frac{r_0+r_0'}{2} .
\end{eqnarray*}
(ii) For a general $x\in M(A,A')$ the existence of a corresponding $y\in M(B,B')$ with $d(x,y)\le \frac{r_0+r_0'}{2}$
just follows by induction and the fact that the convex midpoint map is continuous (see Lemma \ref{lemma-necessary}).
\hfill $\Box$ \\

\begin{proposition} \label{prop-midpoint-2}
Let $(X,d)$ be a metric space and $m:X\times X \longrightarrow X$ be a convex midpoint map for $(X,d)$. Then the map
$\frak{M}:\mathcal{C}_m\times \mathcal{C}_m \longrightarrow \mathcal{C}_m$ defined via
\begin{equation} \label{eqn-prop-midpoint2-proof-1}
\frak{M}(A,B) \; := \; N_{\frac{d_H(A,B)}{2}}A \; \cap \;  N_{\frac{d_H(A,B)}{2}}B \hspace{1cm} \forall A,B\in \mathcal{C}_m
\end{equation}
is a midpoint map for $(\mathcal{C}_m,d_H)$.
\end{proposition}
{\bf Proof of Proposition \ref{prop-midpoint-2}:} From equality (\ref{eqn-prop-midpoint2-proof-1}) it follows that
\begin{equation} \label{eqn-prop-midpoint2-proof-2}
\mathfrak{M}(A,B) \; \subset \; N_{\frac{d_H(A,B)}{2}} (A) \hspace{0.5cm} \wedge \hspace{0.5cm}
\mathfrak{M}(A,B) \; \subset \; N_{\frac{d_H(A,B)}{2}} (B) . 
\end{equation}
With $M$ as in Proposition \ref{prop-midpoint-1} it holds
\begin{displaymath}
d_H\Big( A,M(A,B)\Big) \; = \; \frac{d_H(A,B)}{2} \; = \; d_H\Big( B,M(A,B)\Big) .
\end{displaymath}
Thus we find
\begin{equation} \label{eqn-prop-midpoint2-proof-3}
M(A,B) \; \subset \; N_{\frac{d_H(A,B)}{2}} (A) \hspace{0.5cm} \wedge \hspace{0.5cm}
M(A,B) \; \subset \; N_{\frac{d_H(A,B)}{2}} (B) 
\end{equation}
as well as 
\begin{equation} \label{eqn-prop-midpoint2-proof-4}
A\subset N_{\frac{d_H(A,B)}{2}}\Big( M(A,B) \Big) \hspace{0.5cm} \wedge \hspace{0.5cm}
B\subset N_{\frac{d_H(A,B)}{2}}\Big( M(A,B) \Big) .
\end{equation}
Now (\ref{eqn-prop-midpoint2-proof-1}) and (\ref{eqn-prop-midpoint2-proof-3}) yield
\begin{displaymath}
M(A,B) \; \subset \; \mathfrak{M}(A,B) .
\end{displaymath}
This together with (\ref{eqn-prop-midpoint2-proof-4}) implies
\begin{displaymath}
A\subset N_{\frac{d_H(A,B)}{2}}\Big( \mathfrak{M}(A,B) \Big) \hspace{0.5cm} \wedge \hspace{0.5cm}
B\subset N_{\frac{d_H(A,B)}{2}}\Big( \mathfrak{M}(A,B) \Big) , 
\end{displaymath}
which, combined with (\ref{eqn-prop-midpoint2-proof-2}) yields
\begin{displaymath}
d_H\Big( A,\mathfrak{M}(A,B)\Big) , d_H\Big( B,\mathfrak{M}(A,B)\Big) \; \le \; \frac{d_H(A,B)}{2},
\end{displaymath}
such that the triangle inequality for $d_H$ implies
\begin{displaymath}
d_H\Big( A,\mathfrak{M}(A,B)\Big) \; = \; \frac{d_H(A,B)}{2} \; = \; d_H\Big( B,\mathfrak{M}(A,B)\Big) .
\end{displaymath}
Finally note that the facts that $(X,d)$ is proper and $m$ is convex yield $\mathfrak{M}(A,B)\in \mathcal{C}_m$ .
\hfill $\Box$ \\

It is easy to see that, in contrast to $M$, $\frak{M}$ is not convex in general. \\
Just along the lines of the proof of the Proposition \ref{prop-midpoint-2} one also
achieves the
\begin{proposition} \label{prop-midpoint-3}
Let $(X,d)$ be a proper metric space and $m:X\times X\longrightarrow X$ be a midpoint map for $(X,d)$. Then the map
$\tilde{\mathfrak{M}}:\mathfrak{C}\times \mathfrak{C}\longrightarrow \mathfrak{C}$ defined via
\begin{displaymath}
\begin{array}{rcllc}
\tilde{\mathfrak{M}}(A,B) & := &  N_{\frac{d_H(A,B)}{2}}A \; \cap \;  N_{\frac{d_H(A,B)}{2}}B & \forall A,B\in \mathfrak{C} &
\end{array}
\end{displaymath}
is a  midpoint map for $(\mathfrak{C}(X,d),d_H)$. 
\end{proposition}

Note that Theorem \ref{theo-main} implies that the class of distinguished geodesics in \linebreak
$(\mathcal{C}(X,d),d_H)$ determined via the midpoint map $M$ 
is invariant under any isometry of $(\mathcal{C}(X,d),d_H)$ onto itself. \\
This can be interpreted as a Mazur-Ulam type statement for these metric spaces. Recall that the famous Mazur-Ulam Theorem
(see \cite{mu}) claims that the surjective isomeries from a normed vector space onto itself are linear up to translations, i.e.
that they map straight lines onto straight lines, thus leaving invariant the certain class of geodesics determined by the convex midpoint 
map as given in Example \ref{example-midpoint}. (For an astonishingly nice and simple proof of the Mazur-Ulam Theorem also see \cite{v}).


\section{The proof of Theorem \ref{theo-main}}
\label{sec-proof-1}

In this section we prove Theorem \ref{theo-main}. The strategy of this proof is clearly the same as those given in \cite{sch} and
\cite{g} for the Euclidean case: First we establish that images of points are points, i.e. $i\in Isom(X,d)$ given via
$i(p):=I(p)$ is well defined. Then we prove that the isometry $J\in Isom(\mathcal{C}(X,d),d_H)$ given via
$J(C):=(i^{-1}\circ I)(C)$ for all $C \in \mathcal{C}(X,d)$ is the identity.

\begin{lemma} \label{lemma-proof-1}
Let $(X,d)$ be a geodesic metric space such that geodesics do not split, $p\in X$ and $A,B\in \mathfrak{C}(X,d)$ such that 
\begin{displaymath}
d_H(p,A) \; = \; \frac{1}{2} d_H(A,B) \; = \; d_H(p,B).
\end{displaymath}
Then $\min \{ \# A,\# B\} =1$.
\end{lemma}
{\bf Proof of Lemma \ref{lemma-proof-1}:} Let $h:=d_H(A,p)=d_H(B,p)$. Since $A$ and $B$ are compact, there exist $a\in A$, $b\in B$
such that $d(a,b)=2h$ as well as $d(a,b')\ge 2h$ for all $b'\in B$ or $d(a',b)\ge 2h$ for all $a'\in A$. Let us without loss
of generality assume that $d(a,b')\ge 2h$ for all $b'\in B$. \\
Then $\# B=1$, since for $a\in S_h(p)$ there exists a unique $b\in B_h(p)$ with $d(a,b)\ge 2h$ and thus $B=\{ b\}$.
\hfill $\Box$

\begin{lemma} \label{lemma-proof-2}
Let $(X,d)$ be a uniquely geodesic, geodesically complete metric space such that geodesics do not split and $I\in Isom(\mathcal{C}(X,d),d_H)$ such that
there exists $p\in X$ such that $\# I(p)=1$. Then $\# I(q)=1$ for all $q\in X$. 
\end{lemma}
{\bf Proof of Lemma \ref{lemma-proof-2}:} Let $q\in X$ and choose $\tilde{p}\in X$ such that
\begin{displaymath}
d(\tilde{p},p) \; = \; \frac{1}{2} d(\tilde{p},q) \; = \; d(p,q) .
\end{displaymath}
Then $p$ is the unique midpoint of $\tilde{p}$ and $q$ in $(\mathcal{C}(X,d),d_H)$ and, since $I\in Isom(\mathcal{C}(X,d),d_H)$,
$I(p)$ clearly is the unique midpoint of $I(\tilde{p})$ and $I(q)$ in $(\mathcal{C}(X,d),d_H)$. Therefore $M(I(\tilde{p}),I(q))=I(p)\in X$
with $M$ as defined in Proposition \ref{prop-midpoint-1}, from which follows that $\# I(q)=1$.
\hfill $\Box$

\begin{lemma} \label{lemma-proof-3}
Let $(X,d)$ be a proper, uniquely geodesic, geodesically complete metric space such that geodesics do not split and the midpoint map
is convex. Let further $I\in Isom(\mathcal{C}(X,d),d_H)$, then $\# I(p)=1$ for all $p\in X$.
\end{lemma}
{\bf Proof of Lemma \ref{lemma-proof-3}:} Suppose there exits $A\subset \mathcal{C}(X,d)$ with $\# A>1$ and $\# I(A)=1$, i.e. $I(A) \in X$.
Let $r:=diam A\neq 0$ and $q\in A$ such that there exists $\tilde{q}\in A$ with $d(q,\tilde{q})=r$. \\
For each $x\in I(q)$ we choose the $y(x)$ such that 
\begin{displaymath}
d\Big( I(A),x\Big) \; = \; d\Big( x,y(x)\Big) \; = \; \frac{1}{2} d\Big( I(A), y(x)\Big) ,
\end{displaymath}
set $\tilde{Q}=\bigcup\limits_{x\in I(q)} \{ y(x)\}$ and write $Q$ for the closed convex hull of $\tilde{Q}$. It immediately follows
that $d_H(I(A),I(q))=\frac{1}{2}d_H(I(A),Q)$ and thus
\begin{displaymath}
d_H\Big( Q,I(q)\Big) \; \ge \; d_H(Q,p) \; - \; d_H\Big( p,I(q)\Big) \; = \; \frac{1}{2} d_H(Q,p) . 
\end{displaymath}
In order to see that the opposite inequality, $d_H(Q,I(q))\le \frac{1}{2}d_H(Q,p)$, also holds, we have to show that for all
$z\in Q$ there exists $z'\in I(q)$ such that
\begin{displaymath}
d(z,z') \; \le \; \frac{1}{2} d_H(Q,p) \; = \; d_H\Big( p,I(q)\Big) .
\end{displaymath}
This is obviously true for all $z\in \tilde{Q}$. Next let $z\in Q$ be such that there exist $y_1,y_2\in \tilde{Q}$ with
$d(y_1,z)=d(y_2,z)=\frac{1}{2}d(y_1,y_2)$. Then there exist $x_1,x_2,x\in I(q)$ such that $d(x_1,y_1),d(x_2,y_2)\le d_H(I(A),I(q))$
and $d(x,x_1)=d(x,x_2)=\frac{1}{2}d(x_1,x_2)$. Since $m$ is convex, we derive 
\begin{displaymath}
d(x,z) \; \le \; \frac{1}{2} d(x_1,y_1) \; + \; \frac{1}{2} d(x_2,y_2) \; \le  \; d_H\Big( I(A),I(q)\Big) .
\end{displaymath}
The claim for general $z\in Q$ now follows by induction, applying the same argument again and again, the definition of $Q$ and the
fact that $(X,d)$ is complete. \\
Thus we find 
\begin{displaymath}
d_H(A,q) \; = \; d_H\Big( q,I^{-1}(q)\Big) \; = \; \frac{1}{2} d_H\Big( A,I^{-1}(Q)\Big) 
\end{displaymath}
and it follows from Lemma \ref{lemma-proof-1} that $\# I^{-1}(Q)=1$. \\
Let $p:=I^{-1}(Q)\in X$. Then there exists $z\in A$ such that $d(z,q)=d(q,p)=\frac{1}{2}d(z,p)$ and $q$ is the unique midpoint of 
$p$ and $z$ in $(\mathcal{C}(X,d),d_H)$. Therefore $I(q)$ also is the unique midpoint of $I(p)$ and $I(z)$ in $(\mathcal{C}(X,d),d_H)$.
From Proposition \ref{prop-midpoint-1} and Proposition \ref{prop-midpoint-2} it follows that 
\begin{displaymath}
I(q) \; = \; M\Big( I(z),I(p)\Big) \; = \; \mathfrak{M}\Big( I(z),I(p)\Big)
\end{displaymath}
with $M$ and $\mathfrak{M}$ defined as in the Propositions \ref{prop-midpoint-1} and \ref{prop-midpoint-2}. \\
Next we prove \\
(i) $I(A)\in I(z)$: Since $d_H(I(A),I(q))=d_H(A,q)=r$, it follows that there exists $x\in I(q)$ with $d(I(A),x)=r$. From 
$I(q)=\mathfrak{M}(I(z),I(p))$ we deduce that there exists $z'\in I(z)$ with $z'\in B_r(x)$. Since geodesics do not split
we also know that for all $\tilde{z}\in B_r(x)\setminus I(A)$ we have $d(\tilde{z},y(x))<2r$ for $y(x)$ as in the definition of $Q=I(p)$.
Thus $z'=I(A)$ for otherwise $I(q)=M(I(z),I(p))$ yields the existence of an $x'\in I(q)$ with $d(x',I(A))>r$; a contradiction. \\
Now we establish \\
(ii) $I(A)\in I(q)$: Without loss of generality it holds $\# I(q)>1$ for otherwise the claim of the Lemma follows
from Lemma \ref{lemma-proof-2}. Thus, since $I(q)$ is convex and $B_r(I(A))$ is stricly convex, there exists $x'\in I(q)$ such that
$d(I(A),x')=dist(I(A),I(q))=:r-\epsilon <r$. Suppose now $I(A)\notin I(q)$, i.e. $dist(I(A),I(q))>0$ and denote the midpoint of 
two points $a,b\in X$ by $m(a,b)$. Then, since $I(q)=M(I(z),I(p))$ and $I(A)\in I(z)$, $m(I(A),y(x'))\in I(q)$, but  
$d(I(A),m(I(A),y(x')))$ $=r-2\epsilon$, contradicting $dist(I(A), I(q))=r-\epsilon$. \\
(iii) Now $I(A)\in I(q)$ of course implies $I(A)\in I(p)$. On the other hand, since $r=d_H(z,A)=d_H(I(z),I(A))$, there exists
$z_0\in I(z)$ with $d(z_0,I(A))=r$. Now $m(z_0,I(A))\in I(q)$, from which we conclude $z_0\in I(p)$ and thus
$z_0=m(z_0,z_0)\in I(q)$, due to $I(q)=M(I(z),I(q))$. But then it holds $y(z_0)\in I(p)$ and, once again due to $I(q)=M(I(z),I(p))$,
$m(y(z_0),z_0)\in I(q)$. This, however, contradicts $d_H(I(A),I(q))=r$, since $d(y(z_0),I(A))=\frac{3}{2}r$.
\hfill $\Box$ \\

\begin{lemma} \label{lemma-proof-4}
Let $(X,d)$ be as in Theorem \ref{theo-main}, $I\in Isom(\mathcal{C}(X,d),d_H)$,
$i\in Isom(X,d)$ defined via $i(p):=I(p)$ for all $p\in X$ and $J:=i^{-1}\circ I\in Isom(\mathcal{C}(X,d),d_H)$.
Then for all $p\in X$, $r>0$ it holds
\begin{displaymath}
J\Big( B_r(p)\Big) \; = \; B_r(p) .
\end{displaymath}
\end{lemma}
{\bf Proof of Lemma \ref{lemma-proof-4}:} From the definition of $J$ it follows that $J(p)=p$ for all $p\in X$. Thus we find
\begin{displaymath}
r \; = \; d_H\Big( p,B_r(p)\Big) \; = \; d_H\Big( J(p),J(B_r(p))\Big) \; = \; d_H\Big( p,J(B_r(p))\Big) ,
\end{displaymath}
which yields $J(B_r(p))\subset B_r(p)$. In order to prove the claim, we only have to ensure that $S_r(p)\subset J(B_r(p))$.
Under our assumptions, for all $q\in S_r(p)$ there exists a unique $\tilde{q}\in B_r(p)$ such that $d(q,\tilde{q})=2r$. Now it holds
\begin{displaymath}
2r \; = \; d_H\Big( \tilde{q},B_r(p)\Big) \; = \; d_H\Big( \tilde{q}, J(B_r(p))\Big)
\end{displaymath}
and thus $q\in J(B_r(p))$.
\hfill $\Box$ \\

Now we are ready to provide the \\
{\bf Proof of Theorem \ref{theo-main}:} Let $I\in Isom(\mathcal{C}(X,d),d_H)$ and $i\in Isom(X,d)$ as well as 
$J\in Isom(\mathcal{C}(X,d),d_H)$ be defined as in Lemma \ref{lemma-proof-4}. All we have to prove is that $J(C)=C$ for all $C\in \mathcal{C}(X,d)$. \\
Suppose that there exists $p\in C\setminus J(C)$. Since $J(C)$ is compact, there exists a $q\in J(C)$
such that $dist(p,J(C))=d(p,q)$. Now let $n$ be such that $C,J(C)\subset B_{(2^n-1)d(p,q)}(p)$. Then
\begin{equation} \label{eqn-theo-main-proof}
d_H\Big( J(C),B_{(2^n-1)d(p,q)}(p)\Big) \; = \; d_H\Big( C,B_{(2^n-1)d(p,q)}(p)\Big) \; \le \; (2^n-1)d(p,q).
\end{equation}
Let $p_n\in S_{(2^n-1)d(p,q)}(p)$ be such that $d(p_n,q)=d(p_n,p)+d(p,q)=2^nd(p,q)$. From Lemma \ref{lemma-m-conv-dist}
we know that $dist(J(C),\cdot )$ is $m$-convex. Thus it follows 
\begin{displaymath}
dist \Big( J(C),p_n\Big) \; = \; 2^n d(p,q),
\end{displaymath}
contradicting inequality (\ref{eqn-theo-main-proof}). This proves $C\subset J(C)$. Of course, the same argument with $J$ replaced
by $J^{-1}$ yields $J(C)\subset C$ and therefore $J(C)=C$.
\hfill $\Box$


\section{The proof of Theorem \ref{theo-main-2}}

\label{sec-proof-2}

In this section we prove Theorem \ref{theo-main-2}. This proof is even stronger modeled on the classical one dealing with
the Euclidean space. In fact, once Lemma \ref{lemma-proof2-3} is established in our more general setting, Gruber's original proof
essentially also works in this setting (see \cite{g2}).

\begin{lemma} \label{lemma-proof2-1} (see $(2)$ in \cite{g})
Let $(X,d)$ be a proper, uniquely geodesic metric space, 
$I\in Isom(\mathfrak{C}(X,d),d_H)$ and $p,q\in X$, $p\neq q$. Then
\begin{displaymath}
I(p) \; \subset \; \partial N_{d(p,q)} \Big( I(q)\Big) .
\end{displaymath}
\end{lemma}
{\bf Proof of Lemma \ref{lemma-proof2-1}:} Suppose to the contrary $I(p)\not\subset \partial N_{d(p,q)}(I(q))$. We have 
$d_H(I(p),I(q))=d(p,q)$ and the compactness of $I(p)$ and $I(q)$ yield $I(p)\subset N_{d(p,q)}(I(q))$. Thus for 
$\tilde{\mathfrak{M}}(I(p),I(q))$ as in Proposition \ref{prop-midpoint-3} we deduce \linebreak
$\tilde{\mathfrak{M}}(I(p),I(q))^{\circ}\neq \emptyset$. \\
From $\tilde{\mathfrak{M}}(I(p),I(q))$ we remove a non empty open subset contained in \linebreak
$\tilde{\mathfrak{M}}(I(p),I(q))^{\circ}$
of diameter $<\frac{d(p,q)}{2}$ obtaining a set $D$. Now it is easy to see that $D\neq \tilde{\mathfrak{M}}(I(p),I(q))$
also is a midpoint of $I(p)$ and $I(q)$ in $(\mathfrak{C}(X,d),d_H)$, contradicting the uniqueness of the midpoint of 
$I(p)$ and $I(q)$ in $(\mathfrak{C}(X,d),d_H)$.
\hfill $\Box$ \\

\begin{lemma} \label{lemma-proof2-2}
Let $(X,d)$ be a proper, uniquely geodesic, geodesically complete metric space such that geodesics do not split 
and $I\in Isom(\mathfrak{C}(X,d),d_H)$.
Then $\# I(p)=1$ for all $p\in X$.
\end{lemma}
{\bf Proof of Lemma \ref{lemma-proof2-2}:} Suppose there exists $A\in \mathfrak{C}(X,d)$ with $\# A >1$ and $\# I(A)=1$. 
Then with the notation as in the
proof of Lemma \ref{lemma-proof-3} we find that $I(q)$ is a midpoint of $I(A)$ and $\tilde{Q}$ and there exists $p\in X$ such that 
$I(p)=\tilde{Q}$. Lemma \ref{lemma-proof2-1} applied to $z$ and $p$ as well as to $z$ and $q$ yields $I(A)\in I(z)$, from which together
with Lemma \ref{lemma-proof2-1} follows $I(q)\in S_r(I(A))$. The same argument, of course,
yields $I(z)\in S_r(I(A))$, which clearly contradicts $I(A)\in I(z)$.
\hfill $\Box$ \\

\begin{lemma} \label{lemma-proof2-3}
Let $S\subset S_r(p)$ with $\# S<\infty$. Then, with $J$ defined as in Lemma \ref{lemma-proof-4} it holds $J(S)=S$.
\end{lemma}
{\bf Proof of Lemma \ref{lemma-proof2-3}:} Since $d_H(p,S)=d_H(J(p),J(S))=d_H(p,J(S))$, we find on the one hand
\begin{equation} \label{eqn-lemma2-4-1}
J(S) \; \subset \; B_r(p) .
\end{equation}
On the other hand, it holds 
\begin{equation} \label{eqn-lemma2-4-2}
J(S) \; \cap \; S_r(p) \; = \; S .
\end{equation}
In order to see this, let $q\in S$. Then there exists a unique $\tilde{q}\in B_{r}(p)$ with $d(q,\tilde{q})=2r$. Moreover,
$d(\tilde{q},q')<2r$ for all $q'\in B_r(p)$, $q'\neq q$. But $2r=d_H(\tilde{q},S)=d_H(\tilde{q},J(S))$, hence the inclusion
(\ref{eqn-lemma2-4-1}) implies $q\in J(S)$, which yields $S\subset J(S)$. The opposite inclusion just follows by an
analogue argument interchanging the roles of $S$ and $J(S)$. \\
Furthermore, the same argument yields
\begin{equation} \label{eqn-lemma2-4-3}
S_R(p) \; \subset \; J\Big( B_R(p)\Big) \; \subset \; B_R(p) \hspace{1cm} \forall R\ge 0 .
\end{equation}
From (\ref{eqn-lemma2-4-1}) and (\ref{eqn-lemma2-4-2}) the claim obviously follows, once we establish that
\begin{displaymath}
J(S) \; \cap \; B_r^{\circ}(p) \; = \; \emptyset .
\end{displaymath}
Suppose to the contrary that there exists $q\in J(S)\cap B_r^{\circ}(p)$ and let $\mu := \frac{r-d(p,q)}{2}>0$. Since
$\mu \le \frac{r}{2}$, we find $N_{r-\mu}(S_{\mu}(p))=B_r(p)$. From this and the inclusion (\ref{eqn-lemma2-4-1}) it follows
\begin{equation} \label{eqn-lemma2-4-4}
N_{r-\mu}\Big( S_{\mu}(p)\Big) \; = \; B_r(p) \; \supset \; J(S) ,
\end{equation}
while from $0\le d(p,q)=r-2\mu$ and $q\in J(S)$ we deduce
\begin{equation} \label{eqn-lemma2-4-5}
N_{r-\mu} \Big( J(S)\Big) \; \supset \; B_{r-\mu}(q) \; \supset \; B_{\mu}(p) .
\end{equation}
Now (\ref{eqn-lemma2-4-3}), (\ref{eqn-lemma2-4-4}) and (\ref{eqn-lemma2-4-5}) imply
\begin{displaymath}
d_H\Big( J(S), J(B_{\mu}(p))\Big) \; \le \; r-\mu ,
\end{displaymath}
which contradicts $d_H(S,B_{\mu}(p))\ge r$, since $J$ is an isometry of $(\mathfrak{C}(X,d),d_H)$.
\hfill $\Box$ \\

{\bf Proof of Theorem \ref{theo-main-2}:} We only have to show that for $J$ as in Lemma \ref{lemma-proof2-3} it holds
$J(A)=A$ for all $A\in \mathfrak{C}(X,d)$. \\
Suppose there exists $p\in A\setminus J(A)$. Since $J(A)$ is compact, we have
\begin{displaymath}
\mu \; := \; \frac{1}{2} \inf \Big\{ d(p,q) \; \Big| \; q\in J(A)\Big\} \; > \; 0 .
\end{displaymath}
Thus $U:=\bigcup\limits_{q\in J(A)} B^{\circ}_{d(p,q)-\mu}(q)$ is an open covering of $J(A)$. Since $J(A)$ is compact, there exists
a finite subcovering of $U$, say
\begin{displaymath}
\bigcup\limits_{k=1,...,n} \; B^{\circ}_{d(p,q_k)-\mu} (q_k) \; \supset \; J(A) .
\end{displaymath}
Let $y_1,...,y_k$ be such that $d(y_k,p)=d(p,q_k)+d(q_k,y_k)$, $k=1,...,n$ and 
$d(y_1,p)=d(y_2,p)=...=d(y_n,p)=:\lambda$. Then $\bigcup\limits_{k=1,...,n}B^{\circ}_{\lambda - \mu}(y_k)$ is an open covering of $J(A)$. \\
We set $S:=\{ y_1,...,y_k\} \subset S_{\lambda}(p)$ and obtain, on the one hand
\begin{equation} \label{eqn-proof2-1}
N_{\lambda -\mu}S \; \supset \; J(A) .
\end{equation}
On the other hand it holds $d(q_k,p)\ge 2\mu$ and thus $d(y_k,q_k)=d(y_k,p)-d(p,q_k)\le \lambda -2\mu$, which yields
\begin{equation} \label{eqn-proof2-2}
S \; \subset \; N_{\lambda -2\mu} \Big( J(A)\Big) .
\end{equation}
From Lemma \ref{lemma-proof2-3} we know that $S=J(S)$, which together with the inclusions (\ref{eqn-proof2-1}) and (\ref{eqn-proof2-2})
yields 
\begin{equation} \label{eqn-proof2-3}
d_H\Big( J(S),J(A)\Big) \; \le \; \lambda \; - \; \mu .
\end{equation}
But, since $p\in A$ and $S\subset S_{\lambda}(p)$, we also have $d_H(S,A)\ge \lambda$, contradicting inequality
(\ref{eqn-proof2-3}), due to the fact that $J$ is an isometry of $(\mathfrak{C}(X,d),d_H)$. Hence $A\subset J(A)$ and the same argument replacing 
$J$ through $J^{-1}$ yields $J(A)\subset A$, hence $A=J(A)$ and thus the claim.
\hfill $\Box$ \\


{\footnotesize UNIVERSITY OF MICHIGAN, DEPARTMENT OF MATHEMATICS, 525 E. UNIVERSITY AVENUE, EAST HALL, ANN ARBOR, MI 48109, USA \\
E-mail address: $\;\;\;\;\;$ foertsch@umich.edu}


\end{document}